\documentclass[a4paper,12pt]{article}
\usepackage{amsfonts}
\usepackage{amsmath}
\usepackage{amssymb}
\usepackage{amsthm}
\title{Computation of an Integral Basis of Quartic Number Fields}
\author{Lhoussain El Fadil}
\date{}
\newtheorem{teor}{Theorem}[section]

\newtheorem{lem}[teor]{Lemma}

\newtheorem{rems}[teor]{Remarks}
\newtheorem{prop}[teor]{Proposition}

\newcommand{\z}{\mathbb Z}
\newcommand{\Q}{{\mathbb Q}}

\def\op{\operatorname}

\def\al{\alpha}

\def\be{\bigskip}

\def\dg{\op{deg}}

\def\diso{\lower.4ex\hbox{$\downarrow$}\raise.4ex\hbox{\mbox{\scriptsize $\wr$}}}

\def\fph{\mathbb{F}_{\ph}}

\def\al{\alpha}

\def\ind{\op{ind}}
\def\iso{\,\lower .6ex\hbox{$\stackrel{\lra}{\mbox{\tiny $\sim\,$}}$}\,}

\def\la{\lambda}

\def\lg{l\raise.6ex\hbox to.2em{\hss.\hss}l}
\def\lra{\longrightarrow}

\def\nph#1{N_{\ph}(#1)}
\def\npp#1{N_{\ph}^+(#1)}

\def\orb{\hbox to  .3em{$\backslash$}\backslash}

\def\ph{\phi}

\def\rd{\op{red}}

\def\t{\theta}

\newcounter{cs}
\stepcounter{cs}
\newcommand{\casos}{\begin{itemize}}
\newcommand{\fcasos}{\end{itemize}\setcounter{cs}{1}}

\newfont{\tit}{cmr12 scaled \magstep3}
\begin{document}
\maketitle

\begin{abstract}
In this paper,  based on  techniques of Newton polygons, a result  which allows the computation of a $p$ integral basis of every quartic
number field is given.  For each prime integer $p$, this result allows  to compute a $p$-integral basis of
 a quartic number field $K$ defined by an irreducible polynomial $P(X)=X^4 + aX+b\in
{\z}[X]$ in  methodical and complete generality.
\end{abstract}
\section*{Introduction}
{ Let $K$ be a quartic number field  defined by an irreducible
polynomial $P(X)=X^4 +mX^2+ aX+b\in {\z}[X]$, $\al$ a complex root of
$P$, $\z_K$ the ring of integers of $K$, $d_K$ its discriminant and
$ind(P)=[\z_K:\z[\al]]$ the index of $\z[\al]$ in $\z_K$. It is well
known that: $\triangle=N_{K/{I\!\!\!\!Q}}(P')(\alpha)$ and
$\triangle=(ind(P))^2d_K$, where $\triangle$ is the discriminant of
$P$ and  we can assume that for every prime $p$, $v_p(a)\le
2$ or $v_p(b)\le 3$.\\
 Let $p$ be a prime integer. A $p$-integral basis of $K$ is a
set of integral elements $\{w_1,...,w_4\}$ such that $p$ does not
divide the index $[\z_K:\Lambda]$, where $\Lambda =\sum_{i=1}^{4} \z
w_i $. In that case, we said that $\Lambda$ is a $p$-maximal order
of $K$.  A triangular  $p$-integral basis of $K$ is a $p$-integral
basis of $K$ $\{1,w_2,w_3,w_4\}$ such that $w_1
=\frac{\al+x_1}{p^{r_1}}$, $w_2 =\frac{\al^2+y_2\al+x_2}{p^{r_2}}$
and $w_3 =\frac{\al^3+z_3\al^2+y_3\al+x_3}{p^{r_3}}$. In Theorem
1.2, for
every prime $p$, a triangular  $p$-integral basis of $K$ is given.\\
For every prime $p$ and $(x,m)\in\z^2$, denote
$x_p=\frac{x}{p^{v_p(x)}}$ and $x[m]$: the remainder of the
Euclidean division of $x$ by $m$.\\

In this paper, based on  techniques of Newton polygons, a result which allows the computation of a $p$ integral basis of every quartic
number field is given (Theorem 1.2). Theorem 1.2 allows to find again the results given in \cite{AW} in  methodical and complete generality.
\section*{Newton polygon}
 Let $p$ be a prime integer such that $p^2$ divides $\triangle$
and $\phi(X)$ is an irreducible divisor of $P(X)$ modulo $p$. Set
$m=deg(\phi(X))$ and let
$$
P(X)=a_0(X)\ph(X)^t+a_1(X)\ph(X)^{t-1}+\cdots+a_t(X),
$$ be the $\phi(X)$-adic development of $P(X)$ ( every $a_i(X)\in \z[X]$ and  $\dg a_i(X)<m$). To any coefficient $a_i(X)$ we attach the integer
 $u_i=v_p(a_i(X))$ and the point of the plane $P_i=(i,u_i)$,
 if $u_i<\infty$.

The $\ph$-Newton polygon of $P(X)$ is the lower convex envelope of
the set of points $P_i=(i,u_i)$, $u_i<\infty$, in the cartesian
plane.This (open) polygon is denoted  by $\nph{P}$.

For instance, for a $\ph$-development of degree $7$ with
$u_i=0,\,3,\,0,\,\infty,\,2,\,2,\,4,\,6$ for $i=0,1,\dots,7$, the
polygon is \be
 let $N$ be the $\phi(X)$-Newton polygon of $P(X)$.
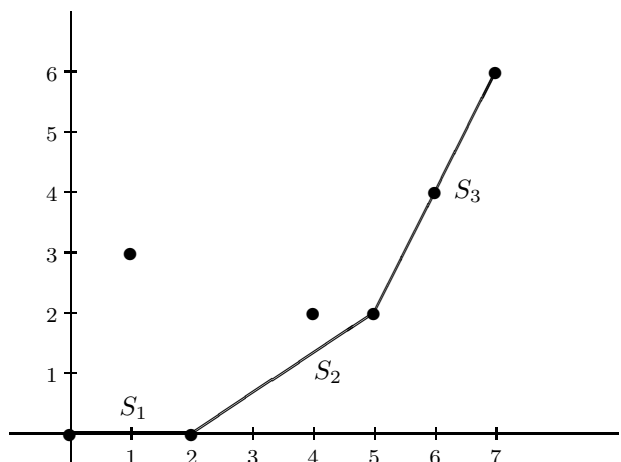
\begin{center}
\setlength{\unitlength}{8.mm}
\begin{picture}(10,6.5)
\put(-.15,-.15){$\bullet$}\put(.85,2.85){$\bullet$}\put(1.85,-.15){$\bullet$}\put(3.85,1.85){$\bullet$}
\put(4.85,1.85){$\bullet$}\put(5.85,3.85){$\bullet$}\put(6.85,5.85){$\bullet$}\put(-1,0){\line(1,0){10}}
\put(0,-1){\line(0,1){8}}\put(1,-.1){\line(0,1){.2}}\put(3,-.1){\line(0,1){.2}}\put(4,-.1){\line(0,1){.2}}
\put(5,-.1){\line(0,1){.2}}\put(6,-.1){\line(0,1){.2}}\put(7,-.1){\line(0,1){.2}}\put(-.1,1){\line(1,0){.2}}
\put(-.1,2){\line(1,0){.2}}\put(-.1,3){\line(1,0){.2}}\put(-.1,4){\line(1,0){.2}}\put(-.1,5){\line(1,0){.2}}
\put(-.1,6){\line(1,0){.2}}\put(2,0){\line(3,2){3}}\put(5,2){\line(1,2){2}}\put(2,0.03){\line(3,2){3}}
\put(5,2.05){\line(1,2){2}}\put(0,0.02){\line(1,0){2}}
\put(.8,.3){\begin{footnotesize}$S_1$\end{footnotesize}}
\put(4.,.9){\begin{footnotesize}$S_2$\end{footnotesize}}
\put(6.3,3.9){\begin{footnotesize}$S_3$\end{footnotesize}}
\put(.9,-.5){\begin{scriptsize}$1$\end{scriptsize}}\put(1.9,-.5){\begin{scriptsize}$2$\end{scriptsize}}
\put(2.9,-.5){\begin{scriptsize}$3$\end{scriptsize}}\put(3.9,-.5){\begin{scriptsize}$4$\end{scriptsize}}
\put(4.9,-.5){\begin{scriptsize}$5$\end{scriptsize}}\put(5.9,-.5){\begin{scriptsize}$6$\end{scriptsize}}
\put(6.9,-.5){\begin{scriptsize}$7$\end{scriptsize}}\put(-.4,.85){\begin{scriptsize}$1$\end{scriptsize}}
\put(-.4,1.85){\begin{scriptsize}$2$\end{scriptsize}}\put(-.4,2.85){\begin{scriptsize}$3$\end{scriptsize}}
\put(-.4,3.85){\begin{scriptsize}$4$\end{scriptsize}}\put(-.4,4.85){\begin{scriptsize}$5$\end{scriptsize}}
\put(-.4,5.85){\begin{scriptsize}$6$\end{scriptsize}}
\end{picture}
\end{center}

The \emph{length} $\ell(\nph{P})$ and the \emph{height}
$h(\nph{P})$ of the polygon are the respective lengths of the
projection to the horizontal and vertical axis. Clearly, $ \dg
P(X)=m\ell(\nph{P})+\dg a_0(X),$ where $m=deg \phi$.\\
 The
$\ph$-Newton polygon is the union of different adjacent
\emph{sides} $S_1,\dots,S_t$ with increasing slope
$\la_1<\la_2<\cdots<\la_t$. We shall write
$\nph{P}=S_1+\cdots+S_t$. The points joining two different sides
are called the \emph{vertexs} of the polygon. The polygon
determined by the sides of positive slopes of $\nph{P}$ is called
the \emph{principal $\ph$-polygon} of $P(X)$ and  denoted by
 $\npp{P}$. The length and the height of $\npp{P}$ are the respective lengths of the
  projection to the horizontal and vertical axis.

For instance, the polygon of the figure has three sides
$S_1,S_2,S_3$ with slopes $0<2/3<2$ and $\npp{P}=S_2+S_3$. For
every side $S$  of the principal part $\npp{P}$, the \emph{length}
$\ell(S)$ and the \emph{height} $h(S)$, of $S$, are the respective
lengths of the projection to the horizontal and vertical axis. The
\emph{slope} of $S$ is the quotient $h(S)/\ell(S)$. The positive
integer $d(S):=\gcd(h(S),\ell(S))$ is called the \emph{degree} of
$S$. Denote  $d:=d(S)$ the degree of $S$, $h:=h(S)/d$ and
$e:=\ell(S)/d$ positive coprime integers such that $h/e$ is the
slope of $S$. Let $s=\lfloor \frac{n}{m}\rfloor$ be the integral
part of $\frac{n}{m}$, where $n=deg(P)$ and $m=deg(\phi)$. For
every $1\le j\le s$, let $H_j$ be the length of the projection of
$P_j(j, u_j)$ to the horizontal axis, $h_j$ its integral part and
$t_j=\rd\left(\dfrac{a_j(X)}{p^{h_j}}\right)$, where $red$ is the
canonical map defined on $\z[X]$ by  reduction modulo $p$. If
$P_j\not \in S$, then $t_j=0$ and  if $P_j \in S$, then
$t_j\neq 0$. If $i$ is the abscissa of the initial point of $S$,
let $P_S(Y)$ be the \emph{residual polynomial} attached to $S$:
$$
P_S(Y):=t_iY^d+t_{i+e}Y^{d-1}+\cdots+t_{i+(d-1)e}Y+t_{i+de}\in\fph[Y].
$$
 Let  $ind_{N}(P): =\sum_{j=1}^sh_j$ the number of points with integer
coordinates that lie below the polygon $N$, strictly above the
horizontal axis  and whose abscissas satisfy $1\le j <l-1$, where
$l$ is the length of $N$, $s=\lfloor \frac{n}{m}\rfloor$, $n=deg(P)$ and $m=deg(\phi)$.\\
 Let $\bar P(X)=(\bar\phi_1(X))^{l_1}\bar P_2$ such that  $\bar\phi_1(X)$ is irreducible modulo $p$, $\bar\phi_1(X)$ does not
 divide $P_2(X)$ and $N_1^+=S_1+..+S_s$ the principal part of
 $N_{\phi_1}(P)$. $P$ is said to be $\phi_1$-regular if for every
 $1\le i\le s$, $P_{S_i}(Y)$ is square free.
 $P$ is said to be $p$-regular if $\bar P(X)=\prod_{i=1}^r{\bar \phi_i}(X))^{l_i}$  the factorization of $\bar
P(X)$ modulo $p$ of irreducible polynomials and for every $1\le i\le r$, $P$ is  $\phi_i$-regular. Since the $p$-regularity of $P$ depends on the choice of $\phi_1(X)$ ,...,$\phi_r(X)$. Then
 $P$ is said to be $p$-regular if there exist $\phi_1(X)$ ,...,$\phi_r(X)$ monic polynomial in $\z[X]$ such that  $\bar P(X)=\prod_{i=1}^r{\bar \phi_i}(X))^{l_i}$ is the factorization of $\bar
P(X)$ modulo $p$  and for every $1\le i\le r$, $P$ is  $\phi_i$-regular.\\
Theorem of index:
$v_p(ind(P))\ge \sum_{i=1}^rm_iind_{N_i}(P)$, where
$m_i=deg(\phi_i(X))$ for every $i$. With equality, if $P(X)$  is a
$p$-regular polynomial. (cf. \cite[p 326]{MN}).
\section{$p$-integral basis of  quartic number field defined by $X^4+aX+b$}

In this section, $K=\Q[\al]$, where $\alpha$ is a complex  root of
an irreducible trinomial $P(X)=X^4+aX+b\in \z[X]$ such that for
every prime $p$, $v_p(a)\le 2$ or $v_p(b)\le 3$ .
 \begin{lem}
Let $p$ be a prime integer and
$w=\frac{t\alpha^3+z\alpha^2+y\alpha+x}{p^i} \in K$. Then
$ch_w=X^4+\frac{A_3}{p^{i}}X^3+\frac{A_2}{p^{2i}}X^2+
 \frac{A_1}{p^{3i}}X+\frac{A_0}{p^{4i}}$ is the characteristic polynomial of $l_w$ the endomorphism of $K$ defined by
 $l_w(x)=wx$, where \\ $A_0=x^4 + 3ax^2yz + 2bx^2z^2 - axy^3 - 4bxy^2z - 3ax^3t + by^4 + b^2z^4 +
b^3t^4 + 3a^2x^2t^2 - 3a^2xyzt + a^2xz^3 - 5abxyt^2 + abxz^2t +
4b^2xzt^2 - a^3xt^3 + 4bx^2yt + 3aby^2zt + 2b^2y^2t^2 - abyz^3 -
4b^2yz^2t + a^2byt^3 - ab^2zt^3$,\\
 $A_1=-(4x^3 - 9ax^2t + 4bxz^2 + 8bxyt + 6axyz
+ 6a^2xt^2 - ay^3 - 4by^2z - 3a^2yzt + a^2z^3 - 5abyt^2 + abz^2t +
4b^2zt^2 - a^3t^3)$,\\ $A_2=6x^2 - 9axt + 3ayz + 4byt + 2bz^2 +
3a^2t^2$ and $A_3=-4x+3at$.\\
In particular,  $w$ is integral if and only if for every $1\le j\le
3$, $\frac{A_j}{p^{ji}}\in \z$.
\end{lem}
The following theorem is an improvement and a specialization of
the theorem of index on quartic number fields.
\begin{teor}
Let $P(X)=X^4+mX^3+nX^2+aX+b\in \z[X]$ be an irreducible polynomial
such that for every prime $p$, $v_p(m)=0$ or $v_p(n)\le1$ or
$v_p(a)\le 2$ or $v_p(b)\le 3$. Let  $p$ be a prime integer. If
$P(X)$ is a $p$-regular polynomial, then we have the following:
\begin{enumerate}
\item
 If $\bar P(X)$ is square free, then $(1,\al, \al^2,\al^3)$
 is a $p$-integral basis of $\z_K$.
 \item
 If $\bar P(X)=(\phi(X))^{4}$, where  $deg \phi=1$, then  $(1,\al,
 \frac{\al^2+a_{3}\al}{p^{h_2}},\frac{\al^3+a_{3}\al^2+a_{2}\al}{p^{h_3}})$
 is a $p$-integral basis of $\z_K$, where $P(X)=\sum_{i=0}^4
a_{i}\phi^i$ is the $\phi_1$-adic development of $P(X)$.
 \item
 If $\bar P(X)=(\phi(X))^{3}P_2$, where  $deg \phi=1$ and $\phi(X)$ does not divide $P_2(X)$, then  $(1,\al,
 \frac{\al^2+a_{3}\al}{p^{h_2}},\frac{\al^3+a_{3}\al^2+a_{2}\al}{p^{h_3}})$
 is a $p$-integral basis of $\z_K$, where $P(X)=\sum_{i=0}^4
a_{i}\phi^i$ is the $\phi$-adic development of $P(X)$.
\item
 If $\bar P(X)=(\phi(X))^{2}P_2$, where  $deg \phi_1=1$, $\phi(X)$ does not divide $P_2(X)$ and $P_2(X)$
 is square free, then  $(1,\al,
 \frac{\al^2+a_{3}\al}{p^{h_2}},\frac{\al^3+a_{3}\al^2+a_{2}\al}{p^{h_3}})$
 is a $p$-integral basis of $\z_K$,  where $P(X)=\sum_{i=0}^4
a_{i}\phi^i$ is the $\phi$-adic development of $P(X)$.
 \item
 If $\bar P(X)=(\phi_1(X))^{2}(\phi_2(X))^{2}$, where  $deg \phi_i=1$ and  $\phi_1(X)$ does not divide $\phi_2(X)$, then
  for every $i$, let $P(X)=\sum_{j=0}^4a_{i,j}\phi_i^j$  be the $\phi_i$-adic development of
  $P(X)$,  $w_i\frac{\al^3+a_{i,3}\al^2+a_{i,2}\al}{p^{h^i_3}}$ and $h_3^i\le h_3^j$.  Then:
  \begin{enumerate}
\item If $h_3^i=0$, then $(1,\al,\al^2,w_j)$ is a $p$-integral basis
of $\z_K$.
\item If $h_3^i\ge 1$, then $(1,\al,w_i-p^{h_3^j-h_3^i}w_j,w_j)$ is a $p$-integral basis
of $\z_K$.
\end{enumerate}
 \item
 If $\bar P(X)=(\phi(X))^{2}$, where $\phi(X)$ is irreducible of degree $2$, then  $(1,\al,
 \frac{\phi(\al)}{p^{h}},\frac{\al\phi(\al)}{p^{h}})$
 is a $p$-integral basis of $\z_K$, where $P(X)=\phi^2+A(X)\phi+B(X)$ is the
          $\phi(X)$-adic development of $P(X)$ and $h$ is the little of $(v_p(A(X)), \lfloor\frac{v_p(B(X))}{2}\rfloor)$.
\end{enumerate}
\end{teor}
{\bf Proof}.
\begin{enumerate}
\item
Case 1.  By Dedekind criterion, since $\bar P$ is square free, then
$(1,\al, \al^2,\al^3)$ is a $p$-integral basis of $\z_K$.
 \item
 Cases 2,3,4. By  theorem of index it suffices to show that every $w_i\in\z_K$,
where  $w_2=\frac{\al^2+a_{3}\al}{p^{h_2}}$ and
$w_3=\frac{\al^3+a_{3}\al^2+a_{2}\al}{p^{h_3}}$. Let
$\phi(X)=X-x_0$. By replacing $P(X)$ by $P(X+x_0)$, we can assume
that $x_0=0$, and then  $w_2=\frac{\al^2+m\t}{p^{h_2}}$ and
$w_3=\frac{\al^3+m\al^2+n\al}{p^{h_3}}$. Let
$Ch_{w_2}(X)=X^4+\frac{2n}{p^{h_2}}X^3+\frac{(n^2+2b+am)}{p^{2h_2}}X^2+
\frac{(amn+bm^2-a^2+2bn)}{p^{3h_2}}X+\frac{(bm^2n-abm+b^2)}{p^{4h_2}}$
 and
 $Ch_{w_2}(X)=X^4+\frac{3a}{p^{h_3}}X^3+\frac{(bn+3a^2)}{p^{2h_3}}X^2+
\frac{(a^3+2abn-b^2m)}{p^{3h_3}}X+\frac{(b^3+a^2bn-ab^2m)}{p^{4h_3}}$
 be  the respective  characteristic polynomial of $l_{w_2}$ and $l_{w_3}$, where
 $l_{w}$ is the endomorphism of $K$ defined by the multiplication by $w$. By definition of $h_2$, $p^{h_2}$ divides $n$. Since $N_{\phi}(P)$
   is convex, then $v_p(b)\ge 3h_2$ and $v_p(a)\ge 2h_2$. Thus, $Ch_{w_2}\in\z[X]$, and then  $w_2\in\z$. For $w_3$, by definition of $h_3$,
    $p^{h_3}$ divides $a$. Since $N_{\phi}(P)$ is convex, then $v_p(b)\ge h_3+( h_3- h_1)\ge h_3+( h_3- h_2)$, and then  $v_p(b)+v_p(n)\ge 2h_3$,
     $3v_p(b)\ge 4h_3$ and $2v_p(b)+v_p(m)\ge 3h_3$. Thus,  $Ch_{w_3}\in\z[X]$ and $w_3\in\z$.
\item
Case 5. As in the previous cases, every $w_k\in\z_K$. By Hensel
lemma, let $(P_1,P_2)\in\z_p[X]^2$ such that $P_1P_2=P$ and $\bar
P_k=\phi_k^2$. Since $\bar P_1$  and $\bar P_2$ are cooprime, then
$v_p(ind(P))=v_p(ind(P_1))+v_p(ind(P_2))$ (cf. \cite{MN}). Since $P$ is
$p$-regular, then for every $k$, $v_p(ind(P_k))=h_3^k$. Thus, If
$h_3^i=0$, then $(1,\al,\al^2,w_j)$ is a $p$-integral basis of
$\z_K$. Else, then for every $k$, let $\phi_k=X-x_k$. Then
$a_{k,3}=4x_k+m$. Since $x_1\neq x_2$ modulo $p$, then $a_{1,3}\neq
a_{2,3}$ modulo $p$, and then
$w_i-p^{h_3^j-h_3^i}w_j=\frac{U(\al)}{p^{h_3^i}}$, where
$U(X)\in\z[X]$ of degree $2$ such that the coefficient of $X^2$ is
cooprime to $p$. Finally, $(1,\al,w_i-p^{h_3^j-h_3^i}w_j,w_j)$ is a
$p$-integral basis of $\z_K$.
     \item Case 6.
          By  Theorem of index it suffices to show that every
          $\phi(\al)\in\z_K$. Let $P(X)=\phi^2+A(X)\phi+B(X)$ be the
          $\phi(X)$-adic development of $P(X)$ and
          $k=v_p(\phi(\al))$. Since $P(\al)=0$, then $2v_p(\phi(\al))\ge v_p(B(\al))$, and  then $k\ge h$. Thus, $\frac{\phi(\al)}{p^{h}}\in\z_K$.
\end{enumerate}

 The following Theorem gives us a triangular $p$-integral basis of
$K$, $v_p(\triangle)$ and $v_p(d_K)$  for every prime integer  $p$.
\begin{teor}
Let $p\ge 5 $ be a prime integer. Under the above hypotheses, a
$p$-integral (resp. a $2$-integral, resp. a $3$-integral) basis   of
${Z\!\!\!Z}_K$ is given in table A (resp. table B, B2*, B3* and
B.3.2, resp. table C)
 $$Table A$$
$$\begin{tabular}{|r|c|c|c|c|r|} \hline
case & $v_p(b)$&$v_p(a)$&$v_p(\triangle)$ &p-integral basis&$v_p(d_K)$\\
\hline
A1&3&$\ge 3$&9&$(1,\alpha,\frac{\alpha^2}{p},\frac{\alpha^3}{p^2})$&3\\
\hline
A2&$\ge 3$& 2&8&$(1,\alpha,\frac{\alpha^2}{p},\frac{\alpha^3}{p^2})$&2\\
\hline
A3& $\ge 1$&0&0&$(1,\alpha,{\alpha^2},{\alpha^3})$& 0\\
\hline
 A4& 2& $\ge 2 $ &6&$((1,\alpha,\frac{\alpha^2}{p},\frac{\alpha^3}{p})$& 2\\
\hline
A5&  1& $\ge 1$ &3&$(1,\alpha,\alpha^2,\alpha^3)$&3\\
\hline
A6&  $\ge 2$& 1 &4&$((1,\alpha,{\alpha^2},\frac{\alpha^3}{p})$&2\\
\hline
A7 & 0& $\ge 1 $ &0&$(1,\alpha,\alpha^2,\alpha^3)$&0\\
\hline
 A8 & 0& 0
&?&$(1,\alpha,\alpha^2,\frac{\alpha^3+t\al^2+t^2\al-3t^3}{p^m})$
 & $v_p(\triangle)[2]$\\
 && 
&&
$m=\lfloor \frac{v_p(\triangle}{2}\rfloor$, $3at+4b=0\,[p^{m+1}]$
 & \\
\hline
\end{tabular}$$
 $$Table B$$
$$\begin{tabular}{|r|c|c|c|r|} \hline
case & conditions&$v_2(\triangle)$ &2-integral basis&$v_2(d_K)$\\
\hline
 B1& $v_2(b)\ge 3$, $v_2(a)= 2$  &8&$(1,\alpha,\frac{\alpha^2}{2},\frac{\alpha^3}{2^2})$&2\\
 \hline
 B2& $v_2(b)=3,  v_2(a)\ge 5$&17&$(1,\alpha,\frac{\alpha^2}{2},\frac{\alpha^3}{2^2})$&11\\
\hline
 B3& $v_2(b)=3$, $ v_2(a)=4$ &16&$(1,\alpha,\frac{\alpha^2}{2},\frac{\alpha^3}{2^2})$&10\\
\hline
B4& $v_2(b)= 3$, $ v_2(a)=3$ &12&$(1,\alpha,\frac{\alpha^2}{2},\frac{\alpha^3}{2^2})$& 6\\
\hline
B5&$b=4+16B$, $a=16A, A=B[2]$& 14 &$(1,\alpha,\frac{\alpha^2+2\al+2}{2^2},\frac{\alpha^3+2\alpha^2+(2+4B)\al}{2^3})  $&4\\
\hline
B6&$b=4+16B$, $a=16A, A\neq B[2]$& 14 &$(1,\alpha,\frac{\alpha^2+2\al+2}{2^2},\frac{\alpha^3+2\alpha^2+2\al}{2^2})  $&6\\
\hline
B7&$b=12+16B$, $a=16A$& 14 &$(1,\alpha,\frac{\alpha^2+2}{4},\frac{\alpha^3+2\alpha}{2^2})  $&6\\
\hline
B8&  $v_2(b)=2$, $ v_2(a)=3 $ &12&$(1,\alpha,\frac{\alpha^2}{2},\frac{\alpha^3+2\al}{4})$&6\\
\hline
B9&  $v_2(b)=2$, $ v_2(a)=2 $ &8&$(1,\alpha,\frac{\alpha^2}{2},\frac{\alpha^3}{2})$&4\\
\hline
B10 & $v_2(b)\ge 2$, $ v_2(a)= 1 $ &4&$(1,\alpha,{\alpha^2},\frac{\alpha^3}{2})$&2\\
\hline
B11 & $v_2(b)=1$, $ v_2(a)\ge 3 $ &11&$(1,\alpha,\alpha^2,\alpha^3)$&11\\
\hline
B12 & $v_2(b)=1$, $ v_2(a)=2 $ &8&$(1,\alpha,\alpha^2,\alpha^3)$&8\\
\hline
B13 & $v_2(b)=1$, $ v_2(a)=1 $ &4&$(1,\alpha,\alpha^2,\alpha^3)$&4\\
\hline
B14 &  $ v_2(a)=0 $ &0&$(1,\alpha,\alpha^2,\alpha^3)$&0\\
\hline
B15 &  $ v_2(a)\ge 3, b=1\,[4] $ &8&$(1,\alpha,\alpha^2,\alpha^3)$&8\\
\hline
B16 &  $ v_2(a)\ge 3, b=3\,[8] $ &8&$(1,\alpha,\frac{\alpha^2+1}{2},\frac{\alpha^3-\al^2+\al-1}{2})$&4\\
\hline
B17 &  $ v_2(a)\ge 3, b=7\,[8] $ &8&$(1,\alpha,\frac{\alpha^2+1}{2},\frac{\alpha^3+4\al^2+6\al}{4})$&2\\
\hline
B18& $ v_2(a)=2, b=1\,[4] $ &9&$(1,\alpha,\alpha^2,\alpha^3)$&9\\
\hline
B19 & $ v_2(a)=2, b=7\,[8] $ &10&$(1,\alpha,\frac{\alpha^2+1}{2},\frac{\alpha^3+\al}{2})$&6\\
\hline
B20& $ v_2(a)=1,  b=3\,[4] $ &4&$(1,\alpha,\alpha^2,\alpha^3)$&4\\
\hline
B21 &  $ v_2(a)=1, b=1\,[4] $ &4&$(1,\alpha,\alpha^2,\frac{\alpha^3-\al^2+\al-1}{2})$&2\\
\hline
B2*& $b=3[8]$, $ v_2(a)=2 $ &*&cf table B2*&*\\
\hline
\end{tabular}$$
If $b=3[4]$ and $ v_2(a)=2 $, then let $A=4+a$, $B=1+a+b$. Consider
$F(X)=P(X+1)=X^4+4X^3+6X^2+AX+B$, $d=A_2^2-3B_2$ and $\t=\al-1$. Then
$$
Table\, B2*  $$
$$\begin{tabular}{|r|c|c|r|} \hline
conditions&$v_2(\triangle)$ &2-integral basis&$v_2(d_K)$\\
\hline
  $v_2(B)+1\ge 2v_2(A)$  &$5+2v_2(A)$&$(1,\al,\frac{\t^2}{2},\frac{\t^3+4\t^2+6\t}{2^m})$, $m=v_2(A)$&3\\
  $v_2(A)=m$ &  & &\\
 \hline
   $v_2(B)=2m$, $v_2(A)=m+1$
 &$2m+7$&$(1,\al,\frac{\t^2}{2},\frac{\t^3+4\t^2+6\t}{2^m})$&5\\
 \hline
  $v_2(B)=2m$, $v_2(A)\ge m+2$ &$2m+8$&$(1,\alpha,\frac{\t^2}{2},\frac{\t^3+4\t^2+6\t}{2^m})$ & 6\\
    \hline
 $v_2(B)+1< 2v_2(A)$, $v_2(B)=1[2]$ &*& cf Table B3* &*\\
 \hline
\end{tabular}$$
$B.3*: v_2(B)+1< 2v_2(A),\, v_2(B)=2k+1$. \\ Let  $t\in\z$ such that $v_2(n_2t+A_2)=s$, $s_1=L(\lfloor
\frac{v_2(d)-1}{2}\rfloor, \lfloor \frac{k+2}{2}\rfloor)$ and
$H(X)=F(X+2^kt)=X^4+m_1X^3+n_1X^2+A_1X+B_1$.
$$
Table\, B3*  $$
 $$\begin{tabular}{|c|c|c|r|}
 \hline
conditions& $v_2(\triangle)$ &2-integral basis&$v_2(d_K)$ \\
\hline
     $2v_2(A)> v_2(B)+1$  &2k+9& $(1,\t,\frac{\t^2}{2},\frac{\t^3+4\t^2+6\t}{2^{k+1}})$&5\\
     $v_2(d)=1$  & & &\\
 \hline
$v_2(A)\ge k+3$&2k+9& $s=1$&3\\
  & &$(1,\t,\frac{\t^2}{2},\frac{\t^3+(-3.2^kt+4)\t^2+6\t-6.2^kt}{2^{k+1}})$&\\
 \hline
 $k=1, v_2(d)\ge 2$&11& $(1,\t,\frac{\t^2}{2},\frac{\t^3+6\t}{2^{2}})$&5\\
 \hline
 $v_2(A)= k+2,\, k\ge 2$&2(k+s+1)+7& $s=s_1$&5\\
 $v_2(d)\ge 3$  & &$(1,\t,\frac{\t^2}{2},\frac{\t^3+(4-3.2^kt)\t^2+(6-4.2^{k+1}t)\t-6.2^kt}{2^{k+s_1+2}})$&\\
 \hline
 $v_2(A)= k+2,  v_2(d)= 2$&2k+11& let $s=1$, and go to B4*.&\\
 \hline
  \end{tabular}$$
   $$ Table B4*$$
$$\begin{tabular}{|c|c|r|} \hline
conditions&2-integral basis&$v_2(d_K)$\\
\hline
  $v_2(B_1)= 2(k+2)$ & $(1,\t,\frac{\t^2}{2},\frac{\t^3+(4-3.2^kt)\t^2+6\t-6.2^kt}{2^{k+2}})$&5\\
 \hline
 $ v_2(B_1)\ge 2(k+2)+1$ &$(1,\t,\frac{\t^2}{2},\frac{\t^3+(4-3.2^kt)\t^2+6\t-6.2^kt}{2^{k+3}})$&3\\
 \hline
  \end{tabular}$$
Tables $B2*$, $B3*$ and $B4*$ could be summarized as follows:
  Let $s\in \z$ such that $P(X+s)$ is $2$-regular: $P(X+s)=X^4+4sX^3+6s^2X^2+AX+B$ ($v_2(A)=r$ and $v_2(B)\ge 2r-1$) or ($v_2(A)\ge r+1$ and $v_2(B)= 2r$). A such $s$ could be computed by iterations:  Let $s_0=1$. If $P(X+s_0)$ is not $2$-regular, then replace $s$ by $s_1=s_0+2^k$. Repeat this iteration until to obtain $P(X+s)$ is $2$-regular. Then $(1,\theta, \frac{\theta^2}{2}, \frac{\theta^3+4s\theta^2+6s^2\theta}{2^r})$ is a 2-integral basis of $\z_K$, where $r=\lfloor\frac{v_2(\triangle)-2-v_2(d_K)}{2}\rfloor$ ($r\ge 2$).
$$
Table\, B*$$
$$\begin{tabular}{|r|c|c|r|} \hline
 conditions&$v_2(\triangle)$ &2-integral basis&$v_2(d_K)$\\
\hline
$v_2(B)\ge 2r-1$, $v_2(A)=r\ge 3$&$2r+5$&$(1,\theta, \frac{\theta^2}{2}, \frac{\theta^3+4s\theta^2+6s^2\theta}{2^r})$&3\\\hline
$v_2(B)= 2r$, $v_2(A)= r+1$&$2r+7$&$(1,\theta, \frac{\theta^2}{2}, \frac{\theta^3+4s\theta^2+6s^2\theta}{2^r})$&5\\\hline
$v_2(B)= 2r$, $v_2(A)\ge r+2$&$2r+8$&$(1,\theta, \frac{\theta^2}{2}, \frac{\theta^3+4s\theta^2+6s^2\theta}{2^r})$&6\\\hline
\end{tabular}$$
$$
Table\, C$$
$$\begin{tabular}{|r|c|c|c|r|} \hline
case & conditions&$v_3(\triangle)$ &3-integral basis&$v_p(d_K)$\\
\hline
C1&$v_3(b)\ge 4$, $ v_3(a)= 2$&11&$(1,\alpha,\frac{\alpha^2}{3},\frac{\alpha^3}{3^2})$&5\\
\hline
 C2& $v_3(b)\ge 4$, $v_3(a)= 1$  &7&$(1,\alpha,\alpha^2,\frac{\alpha^3}{3})$&5\\
 \hline
 C3& $v_3(b)\ge 4,  v_3(a)=0, a^2\neq 1[9]$&3&$(1,\alpha,{\alpha^2},{\alpha^3})$&3\\
\hline
 C4& $v_3(b)\ge 4,  v_3(a)=0, a^2=1[9]$&3&$(1,\alpha,{\alpha^2},\frac{\alpha^3-a\al^2+\al}{3})$&1\\
\hline
 C5& $v_3(b)=3$, $ v_3(a)\ge 2 ,$ &9&$(1,\alpha,\frac{\alpha^2}{3},\frac{\alpha^3}{3^2})$&3\\
\hline
C6& $v_3(b)= 3$, $ v_3(a)= 1$ &7&$(1,\alpha,\alpha^2,\frac{\alpha^3}{3})$& 5\\
\hline
C7&$v_3(b)=3$, $a^2=1[9]$&3&$(1,\alpha,{\alpha^2},\frac{\alpha^3-a\al^2+\al}{3})$&1\\
\hline
C8&$v_3(b)=3$, $v_3(a)=0, a^2\neq 1[9]$ &3& $(1,\alpha,{\alpha^2},{\alpha^3})$&3\\
\hline
C9& $v_3(b)=2$, $ v_3(a)\ge 2 ,$ &6&$((1,\alpha,\frac{\alpha^2}{3},\frac{\alpha^3}{3})$&2\\
\hline
C10& $v_3(b)= 2$, $ v_3(a)= 1$ &6&$(1,\alpha,\alpha^2,\frac{\alpha^3}{3})$& 4\\
\hline
C11&$v_3(b)= 2$, $a^2=1[9]$&3&$(1,\alpha,{\alpha^2},\frac{\alpha^3-a\al^2+\al}{3})$&1\\
\hline
C12&$v_3(b)=2$, $v_3(a)=0$, $a^2\neq 1[9]$&3&$(1,\alpha,{\alpha^2},{\alpha^3})$&3\\
\hline
C13& $v_3(b)=1$, $ v_3(a)\ge 1 ,$ &3&$(1,\alpha,{\alpha^2},{\alpha^3})$&3\\
\hline
C14& $b=6 [9]$, $v_3(a)=0$, $ a^2\neq 4[9] ,$ &3&$(1,\alpha,{\alpha^2},{\alpha^3})$&3\\
\hline
C15& $b=6 [9]$, $ a^2=4[9] ,$ &3&$(1,\alpha,{\alpha^2},\frac{\alpha^3-a\al^2+\al}{3})$&1\\
\hline
C16& $b=3 [9]$, $v_3(a)=0$, $ a^2\neq 7[9] ,$ &4&$(1,\alpha,{\alpha^2},{\alpha^3})$&4\\
\hline
C17& $b=3 [9], a^2=7[9]$, $a^4-a^2+b=0[27]$  & $\ge 6$&$(1,\theta,\frac{\theta^2}{3},\frac{\theta^3+4s\theta^2+6s^2\theta}{3^m})$&$v_3(\triangle)\,[2]$\\
&   & &$as=-4b_3\,[3^{v_3(\triangle)+1}]$&\\
&   & & $\theta=\al-s$, $m=\lfloor\frac{v_3(\triangle)}{2}\rfloor$&\\
\hline
C18& $b=3 [9], a^2=7[9]$, $v_3(a^4-a^2+b)=2$  & $5$&$(1,\alpha,{\alpha^2},\frac{\alpha^3-a\al^2+\al}{3})$&$3$\\
\hline
C19& $v_3(b)= 0$  & $0$&$(1,\alpha,{\alpha^2},{\alpha^3})$&$0$\\
\hline
\end{tabular}$$
\end{teor}
 {\bf Proof}. First, $\triangle=2^8b^3-3^3a^4$ and the proof
is based on the Newton polygon. For every prime $p$, let
$u_3=v_p(a)$, $u_4=v_p(b)$, $\bar P(X)$ the reduction of $P$ modulo
$p$, $N$ the $X$-Newton polygon of $P$  and $N^+$ its principal
part.
\begin{enumerate}
 \item  Case1 ($v_p(b)=3$ and $v_p(a)\ge 2$) : A1, B2, B3, B4, C5.\\
  $\bar{P}(X)=X^4$ and $ind_N(P)=3$.\\
 1) If $u_3\ge 3$, then  $N=S$ is one side and $P_S(Y)=Y+\overline{b_p}$ is square free and $(1,\al,\frac{\al^2}{p}, \frac{\al^3}{p^2})$ is a
 $p$-integral basis of $\z_K$.\\
 2) If $u_3=2$, then  $N=S_1+S_2$ with slopes respectively $1/3$ and
  $1$. $P_{S_1}(Y)=Y+\overline{a_p}$ and $P_{S_2}(Y)=\overline{a_p}Y+\overline{b_p}$ are square free.
  Thus,  $(1,\al,\frac{\al^2}{p}, \frac{\al^3}{p^2})$ is a $p$-integral basis of $\z_K$.
  \item  Case2 ($v_p(b)\ge 2$ and $v_3(a)=1$) : A6, B9, C2, C6, C10.\\
     $N=S_1+S_2$ with slopes respectively $1/3$ and
  $1$. $P_{S_1}(Y)$ and $P_{S_2}(Y)$ are of degree 1.
  So, Thus,  $(1,\al,{\al^2}, \frac{\al^3}{p})$ is a $p$-integral basis of $\z_K$.
  \item  Case3 ($v_p(b)\ge 1$ and $v_3(a)=0$)
  : A3, B3, C3, C4, C7, C8, C11, C12, C14, C15, C16, C17,C18.\\
  If $p\neq 3$,  then $\bar{P}(X)$ is square free and then $v_p(ind(P))=0$.\\
For $p=3$, let $F(X)=P(X-a)=X^4-4aX^3+6a^2X^2+AX+B$, where
$A=-a(4a^2-1)$ and $B=(a^4-a^2+b)$. Then $v_3(A)\ge 1$ and
$v_3(B)\ge 1$. It follows that $v_3(ind(P))=0$ if and only if
$v_3(B)=1$, i.e., if ($a^2=1$ modulo $9$ and $v_3(b)=1$) or
($a^2\neq 1$ modulo $9$ and $v_3(b)\ge 2$, then $v_3(ind(P))=0$.
   Else, then $u_1=0$, $u_2=1$, $u_3=v_3(A)\ge 1$ and $u_4=v_3(B)\ge 2$.
  \begin{enumerate}
  \item
   If $u_3=1$ or  $u_4=2$, then $v_3(ind(P))=1$, $\frac{\theta^3-4a\theta^2+6a^2\theta}{3}$ is integral, and then
   $(1,\al,\al^2,\frac{\al^3-a\al^2+\al}{3})$ is a $3$-integral basis
   of $\z_K$: ($b=6$ and $a^2=4$ modulo $9$) or ($b=3$, $a^2=7$ modulo $9$ and $v_3(a^4-a^2+b)=2$).
    \item
     If $u_3\ge 2$ and   $u_4\ge 3$:($b=3$, $a^2=7$ modulo $9$ and $a^4-a^2+b=0$ modulo $27$), then $v_3(\triangle)\ge 6$. Let $m=\lfloor\frac{v_3(\triangle)-2}{2}\rfloor$, $s\in\z$ such that $as=-4b_3\,[3^{v_3(\triangle)+1}]$ and $\theta=\al-s$. Then $K=\Q[\theta]$ and $f(X)=X^4+4sX^3+6s^2X^2+AX+B$ is the minimal polynomial of $\theta$, where $A=4s^3+a$ and $B=s^4+as+b$. It follows that $u_1=0$, $u_2=1$, $u_3=v_2(\triangle)-3=u_4$. Finally,
    $(1,\theta,\frac{\theta^2}{3},\frac{\theta^3+4s\theta^2+6s^2\theta}{3^m})$ is a $3$-integral basis
   of $\z_K$ and $v_3(d_K)=v_3(\triangle)\,[2]$.
 \end{enumerate}
  \item Case 4 ($v_p(b)=2$ and $v_3(a)\ge 2$): A4, C9, B5, B6, B7, B8:\\
   For $p\neq 2$, $N=S$ is one side and  $P_S(Y)=Y^2+\bar{b_p}$ is square free. Hence $v_p(ind(P))=2$\\
   For $p= 2$, $N=S$ is one side and  $P_S(Y)=(Y+1)^2$. Since $P(X)$ is not $2$-regular, we will use a higher order. Let
   $t_1=(X,1/2,Y+1)$ and $\phi_2(X)=X^2+2$ as defined in\cite{GNM}, page
   16 and let $V_2$ be the $2$-adic valuation of $2^d$-order as defined in\cite{GNM}, page
   17.  Then $V_2(\phi_2)=2$,  $V_2(X)=1$ and for every $x\in \z$,
   $V_2(x)=2v_2(x)$.  Let $P(X)=\phi_2^2(X)-4\phi_2(X)+(4AX+4b_2+4)$,  $R_0=V_2(\phi_2^2(X))=4$,
   $R_1=V_2(4\phi_2(X))=6$ and $R_2=V_2(4AX+4b_2+4)$. From \cite[Th 4.18, p:48]{GNM}, it follows that:
   \begin{enumerate}
    \item
    If $v_2(a)=2$, then $R_2=5$ and $(1,\al,\frac{\al^2}{2},\frac{\al^3}{2})$ is a $2$-integral basis of $\z_K$.
    \item
    If  $v_2(a)= 3$, then $R_2= 7$ and $(1,\al,\frac{\al^2}{2},\frac{\al^3+2\al}{4})$ is a $2$-integral basis of $\z_K$.
    \item
    If $b_2+1=0$ modulo $4$  and $v_2(a)\ge 4$, then $R_2\ge 8$ and $(1,\al,\frac{\al^2+2}{4},\frac{\al^3+2\al}{4})$ is a $2$-integral basis of
$\z_K$.
   \item
    If $b_2+1=2$ modulo $4$ and $v_2(a)\ge 4$, then let $b=4+16B,\, a=16A$, $\phi_2(X)=X^2+2X+2t$ and $P(X)=\phi_2^2(X)-4(X+t-1)\phi_2(X)+8(t-1+2A)X
    +4(1+t^2+4B-2t)$. Then $R_0= 4$ $R_1= 7$. Since $1+t^2+4B-2t=2(1-t+2B)$ modulo $8$, let $t\in \z$ such that
$1-t+2B=0$ modulo $4$. It follows that:
     if $B=A$ modulo $2$, then $R_2\ge 10$, $v_2(ind(P))=5$  and $(1,\al,\frac{\al^2+2\al+2}{4},\frac{\al^3+2\al^2+2(1+2B)\al}{2^3})$ is a
      $2$-integral basis of $\z_K$. \\
    If $B\neq A$ modulo $2$, then $R_2=9$, $v_2(ind(P))=4$  and $(1,\al,\frac{\al^2+2\al+2}{4},\frac{\al^3+2\al^2+2\al}{2^2})$ is a
      $2$-integral basis of $\z_K$.
      \end{enumerate}
\item Case5 ($v_p(b)=1$ and $v_3(a)\ge 1$) : A5, C13, B10,B11, B12.\\
 $N=S$ is one side and $P_S(Y)=Y+\overline{b_p}$ is square
free. Hence $v_p(ind(P))=0$.
\item Case6 $v_p(b)\ge 3$ and $v_3(a)= 2$: A2, C1, B1:\\
  $N=S_1+S_2$ with slopes respectively $2/3$ and
$u_4-2$. Since every $P_{S_i}(Y)$ is square free, then
$v_p(ind(P))=3$.
\item Case7 $v_p(b)=0$ and $v_3(a)\ge 1$: A7, B14, B15,..., B19, B2*, B3* and C19.\\
 If $p\neq 2$, then $\bar{P}(X)$ is square free.\\
For $p= 2$,  according to the Dedekind criterion, let
$f(X)=\frac{P(X)-(X+1)^4}{2}=-2X^3-3X^2+\frac{a-4}{2}X+\frac{b-1}{2}$
and $f(-1)=1+\frac{a}{2}+\frac{b-1}{2}$. Thus, \\
If ($v_2(a)\ge 2$ and $b=1$ modulo $4$) or ($v_2(a)=1$ and $b=3$
modulo $4$), then $(1,\al,\al^2,\al^3)$ is a
      $2$-integral basis of $\z_K$. \\
Else, let $P(X+1)=X^4+4bX^3+6X^2+(4+a)X+(1+a+b)$. Then $u_1=2$, $u_2=1$ and\\
\begin{enumerate}
\item
 If $v_2(a)=1$ and $b=1$ modulo $4$, then $u_3=1$, $u_4\ge
2$ and $(1,\al,\al^2,\frac{\al^3-\al^2-\al-1}{2})$ is a
      $2$-integral basis of $\z_K$.
\item  If ($v_2(a)\ge 3$ and $b=3$ modulo
$8$) or ($v_2(a)=2$ and $b=7$ modulo $8$), then $u_3=u_4= 2$ and
$N=S$ is one side such that $P_S(Y)=Y^2+Y+1$ is irreducible. Hence,
 $(1,\al,\frac{\al^2+1}{2},\frac{\al^3-\al^2-\al-1}{2})$ is a
$2$-integral basis of $\z_K$.
\item  If ($v_2(a)\ge 3$ and $b=7$ modulo
$8$, then $u_3=2$, $u_2=1$, $u_3=2$ and $u_4\ge 3$. Thus,
 $(1,\al,\frac{\al^2+1}{2},\frac{\al^3+2\al1}{4})$ is a
$2$-integral basis of $\z_K$.
\item If $v_2(a)=2$ and $b=3\,[4]$, then let
$F(X)=P(X+1)=X^4+4X^3+6X^2+(4+a)X+(1+b+a)=X^4+4X^3+6X^2+AX+B$, where
$A=4+a$, $B=1+b+a$. It follows that:\\
 If   $b=7\,[8]$, then $N_X(F)=S$ is one side and $F_S(y)=y^2+y+1$. Thus, $(1,\al,\frac{\al^2+1}{2},\frac{\al^3+\al^2+\al+1}{2})$ is a
$2$-integral  basis of $\z_K$.\\
If  $b=3\,[8]$, then 
\begin{enumerate}
  \item  If $v_2(B)+1\ge 2v_2(A)$, then $v_2(\triangle)=5+2v_2(A)$, $v_2(ind(P))=1+v_2(A)$,
  $v_2(d_K)=3$ and  $(1,\al,\frac{\al^2}{2},\frac{\al^3+4\al^2+6\al}{2^k})$ is a $2$-integral
  basis of $\z_K$, where $k=v_2(A)$.
 \item  If $v_2(B)=2k$ is even and $v_2(B)+1< 2v_2(A)$, then \\
 if $v_2(A)\ge k+2$, then $v_2(\triangle)=8+v_2(B)$, $v_2(ind(P))=k+1$, $v_2(d_K)=6$ and
$(1,\al,\frac{\al^2}{2},\frac{\al^3+4\al^2+6\al}{2^k})$ is a
$2$-integral  basis of $\z_K$.\\
  If $v_2(A)=k+1$, then $v_2(\triangle)=2k+7$, $v_2(ind(P))=k+1$, $v_2(d_K)=5$ and
  $(1,\al,\frac{\al^2}{2},\frac{\al^3+4\al^2+6\al}{2^k})$ is a $2$-integral
  basis of $\z_K$.
\item  $v_2(B)$ is odd,  $v_2(B)+1 < 2v_2(A)$ and $v_2(B)\ge 5$. By  Theorem of the polygon $F(X)=H(X)G(X)$
in $\z_2[X]$, where $H(X)=X^2+rX+s$, $G(X)=X^2+RX+S$,\\ $\left\{
\begin{array}{cccc}
v_2(s)=1,& v_2(S)=v_2(B)-1, &v_2(r)\ge
1,&v_2(R)>\frac{v_2(B)-1}{2}\\
B=sR+rS,& A=Sr+Rs,& 4=r+R
\end{array}
\right.$. Since   $4=r+R$, $v_2(r)\ge 3$ and $v_2(r)=2$. So,
$v_2(disc(H))=v_2(r^2-4s)=3$, $v_2(Res(H,G))=2v_2(H(\theta))=2$,
where $\theta$ is a root of $G(X)$. As $v_2(r)\ge 1$ and $v_2(s)=1$,
we have $H(X)$ is irreducible in $\z_2[X].$ On the other hand, as
$disc(F)=disc(H)disc(G)(Res(H,G))^2$,
$v_2(disc(F))=v_2(disc(H))+2v_2(Res(H,G))+v_2(disc(G))=7+v_2(disc(G))$.
Thus,
\begin{enumerate}
\item
If $G(X)$ is irreducible in $\z_2[X]$, then from \cite{MN},
$v_2(ind(F))=0+2+v_2(ind(G))$. Let $\t$ be a root of $G(X)$ and
$u=\frac{\t+x}{2^k}\in \Q_2[\t]$. Since the characteristic
polynomial of $u$ is
$Ch_u=X^2-\frac{2x-R}{2^k}X+\frac{(2x-R)^2+(disc(G))}{2^{2k+2}}$,
where $disc(G)=4S-R^2$, then $u$ is integral if and only if
$2^{2k+2}$ divides $disc(G)$. Therefore, $v_2(ind(G))=\lfloor
\frac{v_2(disc(G))}{2}\rfloor -1$. Thus, $v_2(ind(F))=2+\lfloor
\frac{v_2(disc(G))}{2}\rfloor-1$,  and then if $v_2(\triangle)$ is
even, then $v_2(ind(F))= \frac{v_2(disc(F))-6}{2}$ and $v_2(d_K)=6$.
Else, then $v_2(ind(F))= \frac{v_2(disc(F))-5}{2}$ and $v_2(d_K)=5$.
\item
If $G(X)=(X-\theta_1)(X-\theta_2)$ in
$\z_2[X]$, then $v_2(disc(G))=2v_2(\theta_1-\theta_2))$ and
   $v_2(disc(F))=7+2v_2((\theta_1-\theta_2))$. Let $G_1=X-\theta_1$ and $G_2=X-\theta_2$. Then
$v_2(ind(F))=2+v_2(Res(G_1,G_2)=2+v_2((\theta_1-\theta_2))$, and then
$v_2(d_K)=3$.
\end{enumerate}
In these cases,
$(1,\al,\frac{\al^2}{2},\frac{\al^3+z\al^2+y\al^2+x}{2^m}) $ is a
$2$-integral basis of $\z_K$, where $m=\lfloor
\frac{v_2(disc)-v_2(d_K)}{2}\rfloor -1$ , $x$, $y$ and $z$ are
integers such that $\frac{\al^3+z\al^2+y\al^2+x}{2^m}\in\z_K$.\\
Since we can not compute the coefficients of $G(X)$ in $\Q_p$
neither to test if $G(X)$ is irreducible in $\Q_p[X]$,  we must
give, in $C.1$, a method which  allows us to compute  the integers $x$,
$y$ and $z$ independently of the knowledge of the  irreducibility of
$G(X)$.
\end{enumerate}
\end{enumerate}
\item
$C.1$: Let $F(X)=X^4+4X^3+6X^2+AX+B\in\z[X]$, $\theta=\alpha-1$,
$d=A_2^2-3B_2$, where  $v_2(A)\ge k+2$,  $v_2(B)=2k+1$ and $k\ge 1$. Let $t\in \z$ such that
   $v_2(3t+A_2)=s$,
    $H(X)=F(X+2^kt)=X^4+m_1X^3+n_1X^2+A_1X+B_1$. Then   $m_1=4+2^{k+2}t$, $n_1=6+3.2^{k+2}t+3.2^{2k+1}t^2$,
   $A_1=A+3.2^{k+2}t+3.2^{2k+2}t^2+2^{3k+2}t^3$ and   $B_1=B+2^kAt+3.2^{2k+1}t^2+2^{3k+2}t^3+2^{4k}t^4$. Therefore $3B_1=2^{2k+1}(3B-A_2^2)+3.2^{3k}t^3(2^2+2^{k}t)+2^{2k+1+2s}L)$, where $L$  is an odd integer. It follows that:  $v_2(m_1)= 2$,   $v_2(n_1)=1$ and
      \begin{enumerate}
        \item
        If $v_2(d)=1$, then  $v_2(B_1)=2k+2$, $v_2(ind(F))=k+2$, and
    $(1,\t,\frac{\t^2}{2},\frac{\t^3+4\t^2+6\t}{2^{k+1}})$ is  a $2$-integral basis of
    $\z_K$. Moreover, if $v_2(A)=k+2$, then  $v_2(\triangle)=2k+9$ and $v_2(d_K)=5$. if $v_2(A)\ge k+3$, then  $v_2(\triangle)=2k+10$ and $v_2(d_K)=6$.
         \item
         If $v_2(d)\ge 2$ and $v_2(A)>k+2$, then $v_2(A_1)=k+2$ and
          $v_2(B_1)\ge 2k+3$. Hence   $v_2(ind(F))=k+3$,  $v_2(\triangle)=2k+9$ and $v_2(d_K)=3$ and
          $(1,\t,\frac{\t^2}{2},w_3)$ is
   a $2$-integral basis of $\z_K$, where $w_3=\frac{\t^3+(4-3.2^kt)\t^2+6\t-3.2^{k+1}t}{2^{k+2}}$.
   \item
         If  $k=1$, $v_2(d)\ge 2$ and $v_2(A)=3$, then for $s=1$, $v_2(B_1)=4$, $v_2(\triangle)=11$, $v_2(d_K)=5$,
          $(1,\t,\frac{\t^2}{2},\frac{\t^3+4\t^2+6\t}{2^{2}})$ is
   a $2$-integral basis of $\z_K$.
  \item
   If $v_2(A)=k+2$, $k\ge 2$ and  $v_2(d)\ge 3$, then for
   $s=L(\lfloor\frac{v_2(d)-1}{2}\rfloor,\lfloor\frac{k}{2}\rfloor)$,
   $v_2(A_1)= (k+1+s)+1$, $v_2(B_1)\ge 2(k+s+1)+1$. Hence
    $v_2(ind(F))=k+2+s$,  $v_2(\triangle)=2(k+s+1)+7$, $v_2(d_K)=5$
     and $(1,\t,\frac{\t^2}{2},\frac{\t^3+(4-3.2^kt)\t^2+6\t-3.2^{k+1}t}{2^{k+1+s}})$ is
   a $2$-integral basis of $\z_K$.
   \item
   If $v_2(A)=k+2$ and  $v_2(d)= 2$, then $v_2(\triangle)=2k+11$. For
   $s=1$,  $v_2(A_1)= k+3$ and  $3B_1=-2^{2k+1}d+3.2^{3k+2}t^3(1+2^{k-2}t)+2^{2k+3}L)=2^{2k+3}(-d_2+L)+3.2^{3k+2}t^3(1+2^{k-2}t)$, where $L$  is an odd integer. Thus,  $v_2(A_1)= k+3$ and $v_2(B_1)\ge 2(k+2)$. It follows that:
    \begin{enumerate}
    \item
    If $v_2(B_1)= 2(k+2)$, then $v_2(ind(F))=k+3$, $v_2(d_K)=5$  and
   $(1,\t,\frac{\t^2}{2},\frac{\t^3+(4-3.2^kt)\t^2+6\t-3.2^{k+1}t}{2^{k+2}})$ is
   a $2$-integral basis of $\z_K$.
   \item
     If $v_2(B_1)\ge 2(k+2)+1$, then $v_2(ind(F))=k+4$, $v_2(d_K)=3$ and
   $(1,\t,\frac{\t^2}{2},\frac{\t^3+(4-3.2^kt)\t^2+6\t-3.2^{k+1}t}{2^{k+3}})$ is
   a $2$-integral basis of $\z_K$.
    \end{enumerate}
    \end{enumerate}
  \item Case8 $v_p(ab)=0$: A8.\\
   If $p\in\{2,3\}$, then $v_p(\triangle)=0$, $(1,\al,\al^2,\al^3)$ is a $p$-integral basis of $\z_K$  and $v_p(d_K)=0$.\\
Let $p\ge 5$. If $v_p(disc(2^8b^3-3^3a^4))\le 1$, then
$(1,\al,\al^2,\al^3)$ is a $p$-integral basis of $\z_K$  and $v_p(d_K)=0$.\\
Else, since $3a\neq 0$ modulo $p$, let $t\in\z$ such that $3at+4b=0$
modulo $p^s$, where $s=m+1$ and $m=\lfloor
\frac{v_p(\triangle)}{2}\rfloor$ ($3at+4b=p^sL)$. Then
$(3a)^3P'(t)=-\triangle+3.4^3b^2p^sL$ modulo $p^{2s}$. Thus,
$v_p(P'(t))\ge s$. Moreover, $(3a)^4P(t)=b\triangle-p^s\triangle$
modulo $p^{2s}$. Thus, $v_p(P(t))=v_p(\triangle)$.  Let
$P(X+t)=X^4+4tX^3+6t^2X^2+P'(t)X+P(t)$. Since $6t^2\neq 0$ modulo
$p$, then $N=S_0+S_1$ with slopes respectively $0$ and
$\frac{v_p(\triangle)}{2}$ and $P_{S_1}(Y)$ is square free. Hence,
$v_p(ind(P))=\lfloor \frac{v_p(\triangle)}{2}\rfloor$,
$\frac{\t^3+4t\t^2+6t^2\t}{p^m}\in\z_K$, where $\t=\al-t$. Thus,
$(1,\alpha,\alpha^2,\frac{\alpha^3+t\al^2+t^2\al-3t^3}{p^m})$ is a
$p$-integral basis of $z_K$ and $v_p(d_K)=v_p(\triangle)$ modulo
$2$.
\item Case9 : $v_2(a)=0$ or $v_3(b)=0$ (B13,  C19). Since $v_p(\triangle)=0$, $(1,\al,\al^2,\al^3)$ is a $p$-integral basis of $\z_K$  and $v_p(d_K)=0$.
\end{enumerate}
\section{ An integral basis of a quartic number field defined by $X^4+aX+b$}
\begin{rems}
\begin{enumerate}
\item
Let $p$ be a prime integer such that $p^2$ divides $\triangle$. For
every $1 \le i\le 3$, let
$w_{i,p}=\frac{L_i^p(\alpha)}{p^{r_{i,p}}}$, where $L_i^p(X)\in
\z[X]$ is a monic polynomial of degree $i$ such that ${\cal
F}=(1,w_{1,p},w_{2,p},w_{3,p})$ is a triangular $p$-integral basis
of $K$. Ten $r_{1,p}\le r_{2,p}\le r_{3,p}$, $v_p(\triangle)=r_1+
r_2+ r_3$ and $v_p(d_K)=v_p(\triangle)-2(r_1+ r_2+ r_3)$.
\item
 Let $p_1,\cdots,p_r$ be the primes  such that every
$p_i^2$ divides $\triangle$. For every $1 \le i\le 3$, denote
$d_i=\prod_{j=1}^rp_j^{r_{ij}}$, where for every $j$,
$w_{i,j}=\frac{L_i^{p_j}(\alpha)}{p^{r_{ij}}}$ and
$(1,w_{1,j},w_{2,j},w_{3,j})$ is a $p_j$-integral basis of $K$. Then
$1\mid d_1\mid d_2\mid d_3$ are the elementary divisors of
$\z_K/\z[\al]$. In particular, $d_3$ is the conductor of the order
$\z[\al]$ and $ d_1d_2d_3=\mp \ind(P)$.
\item
 We can always assume that a triangular $p$-integral basis has the property:
 if $ r_i=r_{i+1}$, then we can take $w_{i+1}=\al w_i.$
 \end{enumerate}
\end{rems}
One can recover a triangular integral basis from different
triangular $p$-integral basis for all $p$ as follows :
\begin{prop}
Let $p_1$,...,$p_s$ the prime integers such that $p^2$ divides
$\triangle$ and $1,\, d_1,\, d_2$ and $d_3$ the elementary divisors
of the abelian group $\z_K/\z[\al]$. For every $j$, let ${\cal
F}_j=(1,w_{1,j},w_{2,j},w_{3,j})$ be a triangular $p_j$-integral
basis of $K$, i.e., $w_{i,j}=\frac{L_i^j(\alpha)}{p_j^{r_{ij}}}$
such that every $L_i^j(X)$ is a monic polynomial of $\z[X]$ of
degree $i$.
 Then ${\cal B}=(1,w_1,w_2,w_3)$ is a triangular integral basis of $K$, where
every $w_i=\frac{L_i(\alpha)}{d_i}$, $L_i(X)=L_i^{j}(X)$ modulo
$p_j^{r_{ij}}$.
\end{prop}
{\bf Proof} Since $ind(P)=d_1d_2d_3$, we need only to check that
every $w_i\in \z_K$.  Let $1\le i\le 3$. Since for every $i$ the
integers $(\frac{d_i}{p_j^{r_{ij}}})_{1\le j\le s}$ are pairwise
coprime, there exist integers $t_1,...,t_s$ such that
$\sum_{j=1}^st_j\frac{d_i}{p_j^{r_{ij}}}=1$. Hence, $
\frac{L_i(\al)}{d_i}=\sum_{j=1}^st_j\dfrac{L_i(\al)}{p_j^{r_{ij}}}\in
 \z_K,$ because all $\frac{L_i(\al)}{p_j^{r_{ij}}}\in \z_K$.\\

 {\bf Acknowledgments}\\
 This work is supported by the Spanish Ministry of Higher Education (Ref: SB 2006-0128). I would like to thank CRM at Barcelona for their extraordinary hospitality and facilities. I am also thankful to Professor E. Nart for his valuable comments and suggestions.

L HOUSSAIN EL FADIL,  FPO, P.O. Box 638-Ouarzazat\\
\hspace*{3cm} Morocco\\

lhouelfadil@hotmail.com
\end{document}